\documentclass[12pt]{amsart} %\documentclass[11pt]{amsart} %\documentclass[12pt,a4paper]{article}
\usepackage{amsfonts}
\usepackage{amsmath}
\usepackage{amssymb}
\usepackage{amsthm}
\usepackage{graphicx}
\usepackage[T1]{fontenc}

\usepackage{hyperref}

\usepackage{xcolor}

\newtheorem{theorem}{Theorem}
\newtheorem{conjecture}[theorem]{Conjecture}
\newtheorem{lemma}[theorem]{Lemma}
\newtheorem{observation}[theorem]{Observation}
\newtheorem{corollary}[theorem]{Corollary}

\newtheorem{case}{Case}
\newtheorem{claim}{Claim}

\begin{document}

	\title{On list extensions of the majority edge colourings}

	\author{Pawe{\l} P\k{e}ka{\l}a}
	\address[agh]{AGH University of Krakow, al. A. Mickiewicza 30, 30-059 Krakow, Poland}
	\email{\tt ppekala@agh.edu.pl}
	
	\author{Jakub Przyby{\l}o}
	\address[agh]{AGH University of Krakow, al. A. Mickiewicza 30, 30-059 Krakow, Poland}
	\email{\tt jakubprz@agh.edu.pl}

	\begin{abstract}
	We investigate possible list extensions of generalised majority edge colourings of graphs and provide several results concerning these.
	Given a graph $G=(V,E)$, a list assignment $L:E\to 2^C$ and some level of majority tolerance $\alpha\in(0,1)$, an $\alpha$-majority $L$-colouring of $G$ is a colouring $\omega:E\to C$ from the given lists such that for every $v\in V$ and each $c\in C$, the number of edges coloured $c$ which are incident with $v$ does not exceed $\alpha\cdot d(v)$.	
	We present a simple argument implying that for every integer $k\geq 2$, each graph with 
	minimum degree $\delta\geq 2k^2-2k$ admits a $1/k$-majority $L$-colouring from any assignment of lists of size $k+1$. This almost matches the best result in a non-list setting and solves a conjecture posed for the basic majority edge colourings, i.e. for $k=2$, from lists.
	We further discuss restrictions which permit obtaining corresponding results in a more general setting, 
	i.e. for diversified $\alpha=\alpha(c)$ majority tolerances for distinct colours $c\in C$.
        Consider a list assignment $L:E\to 2^C$ with $\sum_{c\in L(e)}\alpha(c)\geq 1+\varepsilon$ for each edge $e$, and suppose that 
        $\alpha(c)\geq a$ for every $c$ or $|L(e)|\leq\ell$ for all edges $e$, where $a\in(0,1)$, $\varepsilon>0$, $\ell\in\mathbb{N}$ are any given constants.
        Then we in particular show that there exists an $\alpha$-majority $L$-colouring of $G$ from any such list assignment, provided that 
         $\delta(G)=\Omega(a^{-1}\varepsilon^{-2}\ln(a\varepsilon)^{-1})$ or $\delta=\Omega(\ell^2\varepsilon^{-2})$, respectively.
	We also strengthen these bounds within a setting where each edge is associated to a list of colours with a fixed vector of majority tolerances, applicable also in a general non-list case.

	\end{abstract}

			\maketitle

	\section{Introduction}

	A \emph{majority colouring} of a graph $G=(V,E)$ is an assignment $\omega:V\to C$ such that at most half the vertices adjacent with any given $v\in V$ are coloured $c(v)$. 
	Research devoted to majority colourings, though under different terminology, dates 
	back to a classic result of Lov\'asz~\cite{Lovasz-majority}, who observed in 1966 that every graph is majority $2$-colourable, or more generally, $1/k$-majority $k$-colourable, where only $1/k$'th fraction of neighbours of every given vertex  may have its colour assigned.
	This concept was further developed in particular  in concern to intriguing 
	problems regarding infinite graphs, see e.g.~\cite{Berger,Bruhn-Diestel,list_maj,Shelah-Milner}.
	The term `\emph{majority colouring}' was later introduced by Kreutzer, Oum, Seymour, van der Zypen and Wood~\cite{digraph1} in the environment of vertex-colourings of \emph{digraphs}, where each vertex was allowed to have at most half the out-neighbours in its own colour. 
	In~\cite{digraph1} it was in particular proven that every digraph admits such a $4$-colouring, and conjectured that $3$ colours should always be enough.
	See~\cite{Szabo-majority,digraph2,MajorityGeneralOur,digraph3,Knox-Samal,digraph1} for various extensions of this concept and a list of corresponding results. 
	
	In~\cite{majority23} a number of similar problems and questions were raised in the context of edge colourings.
	These are also 	somewhat related with other concepts, investigated e.g. in~~\cite{defective,frugal}, 
	and have a slightly different nature than their counterparts in the vertex setting. 
	Consider a graph $G=(V,E)$ and a vertex $v\in V$. By $d_G(v)$ or simply $d(v)$ we denote the degree of $v$ in $G$. Given any $F\subseteq E$, we moreover denote by $F(v)$ the set of all edges belonging to $F$ which are incident with the vertex $v$, and set $d_F(v):=|F(v)|$. Let $\omega:E\to C$ be an edge colouring of $G$. Denote by $E_c$ the set of edges of $G$ coloured $c$. The colouring $c$ is said to be a \emph{majority edge colouring} of $G$ if for every colour $c\in C$, each vertex $v\in V$ is incident with at most $0.5d(v)$ edges coloured $c$, i.e.  $d_{E_c}(v)\leq 0.5d(v)$.
	Note that unlike in the case of colouring vertices, such an edge  colouring does not exist (regardless of the size of $C$) for some graphs, in particular for graphs having a vertex of degree $1$.
	% and some graphs with  the minimum degree $2$ (e.g. odd cycles).
	On the other hand, in~\cite{majority23} it was proven that any graph with minimum degree $\delta\geq 2$ admits a majority edge $4$-colouring, whereas a majority edge $3$-colouring is admissible for any graph with $\delta\geq 4$ (an alternative proof of the latter fact may be found in~\cite{majority_ejc}). Many efforts within research on generalisations of majority edge colourings are thus devoted to establishing which minimum degrees 
	guarantee the existence of such colourings (with a given number of colours). Consider an integer $k\geq 2$. We say that  $\omega:E\to C$ is a \emph{$1/k$-majority edge colouring} if $d_{E_c}(v)\leq (1/k)d(v)$ for every $v\in V$ and $c\in C$. 
	Note such a colouring cannot exist if a graph has a vertex with degree smaller than $k$ or if we strive to use less than $k$ colours. Moreover, it is straightforward to notice that there are graphs with arbitrarily large minimum degrees which do not admit a $1/k$-majority edge $k$-colouring, it suffices if these have a vertex with degree indivisible by $k$. 
	\begin{theorem}[\cite{majority23}]\label{Th1kKalinowski}
	For every integer $k\geq 2$ there exists $\delta^{(k)}=O(k^3\log k)$ such that every graph $G$ with minimum degree $\delta\geq \delta^{(k)}$ is $\frac{1}{k}$-majority edge $(k+1)$-colourable.
	\end{theorem}
	Let us remark that
	the authors of Theorem~\ref{Th1kKalinowski} do not specify the order of $\delta^{(k)}$  within the theorem itself, only mention that it can be 
	%In fact the authors of Theorem~\ref{Th1kKalinowski} do not specify the order of $\delta^{(k)}$ in~\cite{majority23}, but it can be easily 
	derived from the calculations within their probabilistic argument.
	Nevertheless, the lower bound  in Theorem~\ref{Th1kKalinowski} occurred to be rather far from optimal,
in view of the following result.
	\begin{theorem}[\cite{majority_ejc}]\label{thm1-GenMajEdg}
		For every integer $k\geq 2$, if a graph $G$ has minimum degree $\delta\geq 2k^2$, then $G$ is $\frac{1}{k}$-majority edge $(k+1)$-colourable.
	\end{theorem}
	This settled the right order of magnitude for the optimal minimum value of 	$\delta^{(k)}$, introduced in Theorem~\ref{Th1kKalinowski}, as witnesses the following observation.
	\begin{observation}[\cite{majority_ejc}]\label{ExampleGeneral}
		For every $k\geq 2$ there exists a graph $G$ with minimum degree $\delta= k^2-1$ which is not $\frac{1}{k}$-majority edge $(k+1)$-colourable. 
	\end{observation}
Theorem~\ref{thm1-GenMajEdg} was also strengthened for $k\geq 3$ to the following form
(see Theorem 14 in~\cite{majority_ejc} for a slightly stronger result, which implies the one below).
	\begin{corollary}[\cite{majority_ejc}]\label{7_4_Corollary}
		For every integer $k\geq 2$, if a graph $G$ has minimum degree $\delta\geq \frac{7}{4}k^2+\frac{1}{2}k$, then $G$ has a $\frac{1}{k}$-majority edge $(k+1)$-colouring.
	\end{corollary}	
	It was also conjectured that the following holds true. Note that in view of Observation~\ref{ExampleGeneral},
	we cannot expect anything more. 
	\begin{conjecture}[\cite{majority_ejc}]\label{EdgeMajconj}
		For every integer $k\geq 2$, if a graph $G$ has minimum degree $\delta\geq k^2$, then $G$ is $\frac{1}{k}$-majority edge  $(k+1)$-colourable.
	\end{conjecture}
Conjecture~\ref{EdgeMajconj} was confirmed for $k=2$ in~\cite{majority23}, as mentioned above, and for  $k=3,4$ in~\cite{majority_ejc}.	The corresponding problem was also entirely solved (for all $k$) in the case of bipartite graphs~\cite{majority_ejc}.
	
	A list variant of majority edge colourings was first considered  by Kalinowski, Pil\'sniak and Stawiski in~\cite{list_maj}, whose major matter of concern were infinite graphs. We shall focus on the finite setting, though.
	Let $G=(V,E)$ be a  graph. Given a set of colours $C$, we call $L:E\to 2^C$ an $\ell$-list assignment if $|L(e)|=\ell$ for every $e\in E$. An $L$-colouring is an assignment $\omega:E\to C$ such that $\omega(e)\in L(e)$ for every edge $e\in E$.
		\begin{theorem}[\cite{list_maj}]\label{4-list-th}
		Let $G$ be a graph of arbitrary order and without pendant edges. Then there is 
    a majority edge $L$-colouring for any $4$-list assignment $L$.
	\end{theorem}
	It was however conjectured that the following strengthening of Theorem~\ref{4-list-th} should hold.
	\begin{conjecture}\label{list_maj_conj_k2}
		Every graph with minimum degree at least $4$ admits a majority edge colouring from lists of size $3$.
	\end{conjecture}	
Note this would be a direct generalization of the mentioned result from~\cite{majority23} concerning the non-list setting. In the next section we provide an argument confirming Conjecture~\ref{list_maj_conj_k2} (in the finite case), see Theorem~\ref{thm_galvin}.
We actually prove the corresponding result in general for all $k\geq 2$ and graphs
with minimum degree $\delta=\Omega(k^2)$, which settles the order of magnitude in this problem, due to Observation~\ref{ExampleGeneral}.
 It is worth mentioning here that the original probabilistic proof of Theorem~\ref{Th1kKalinowski} 
could be modified towards the discussed list setting.
This would however guarantee a result corresponding to the one presented in the next section, but with a much worse bound: $\delta=\Omega(k^3\log k)$ (with also an inferior, larger multiplicative constant than the one in
Theorem~\ref{Th1kKalinowski}). 
	
	In the third section we observe that our result from Theorem~\ref{thm_galvin} extends to a slightly more general setting, where every vertex may be incident with $\alpha d(v)$ monochromatic edges, where $\alpha$ is any fixed \emph{real} number, called a tolerance. Next, in Section~\ref{Section4}, we discuss reasonable boundaries of possible further extensions of our list setting, in particular towards admitting diversified 
	tolerances for distinct colours around every vertex. The following section is devoted to a 
	resulting most general model admitting `reasonable results'. We exploit the probabilistic method in that and in the following section, which in turn contains an 	improvement of the previous results in the general model with 
an additional restriction, boiling down to unifying multisets of tolerances represented in each list.
This in particular implies certain best known bounds within a non-list setting with diversified tolerances.
We close the paper with a section containing concluding remarks and several conjectures.

	\section{$1/k$-majority list edge colourings}

	In this section we strengthen and extend in Theorem~\ref{thm_galvin} the result from Theorem \ref{thm1-GenMajEdg} to the list setting.
	We provide a straightforward proof, based on the famous Galvin's theorem confirming the List Colouring Conjecture for bipartite graphs. Denote by $\chi'_l(G)$ the list chromatic index of a graph $G$.
	\begin{theorem}[Galvin's theorem, \cite{galvin}]
	For every bipartite graph $G$, $\chi'_l(G)=\chi'(G)=\Delta(G)$.
	\end{theorem}
We shall also exploit an operation of vertex splitting, 
which is a rather standard technique in the environment of majority edge colourings,
cf.~\cite{majority23,majority_ejc} for its previous applications. Let us remark, that unless stated otherwise, by a graph we shall always mean a finite simple graph.
	
	\begin{theorem}\label{thm_galvin}
		For every integer $k\geq 2$, each graph $G$ with minimum degree $\delta \geq 2k^2-2k$ has a $1/k$-majority edge colouring from lists of size $k+1$.
	\end{theorem}
	
%	{\color{magenta} Luzne skojarzenie: Sprawdzic, czy jakies wyniki dotyczace jader nie daja tu czegos, bo byly przydatne w badaniach nad List Colouring Conjecture...}

	\begin{proof}
		Let $G=(V,E)$ be a graph with minimum degree $\delta\geq 2k^2-2k$. 
		We first split every vertex into vertices with degrees at most $2(k+1)$.
		More precisely, 	given any vertex $v$ of $G$, its degree can be written as $d(v)=2(k+1)(s-1)+i$ where $i\in[2k+2]$, i.e. $s=\lceil d(v)/(2k+2)\rceil$. Let us partition the neighbourhood of $v$ in $G$ into $s$ disjoint sets $N_1, \dotsc, N_s$ such that $|N_1|=i$ and $|N_j|=2(k+1)$ for $j\geq 2$. Let $G'$ be a graph such that $V(G') = (V\smallsetminus\{v\}) \cup \{v_1,\dotsc,v_s\}$  and $E(G') = E(G-v)\cup \bigcup_{i=1}^{s} \{uv_i : u\in N_i\}$, where $v_1,\dotsc,v_s$ are new vertices, which can be thought of as copies of $v$ resulting from its splitting. Note that this operation yields a natural bijection between the edges of $G$ and the edges of $G'$. Let $\overline{G}$ be a graph constructed from $G$ by applying the above operation to all vertices of $G$, one after another. If $\overline{G}$ contains a vertex of odd degree, we 
	take an isomorphic copy $\overline{G}'$ of $\overline{G}$	and join
	by an edge every vertex $v$ of odd degree in $\overline{G}$ with its corresponding vertex in $\overline{G}'$.
 (If any extra edges were added in this way, assign to them arbitrary lists of size $k+1$).
We denote the final resulting graph by $G^*$. Note that $\Delta(G^*)=2(k+1)$.
		
		As all vertices of $G^*$ have even degrees, we can find an Euler tour in each component of $G^*$. By giving each of those tours an orientation we obtain a directed graph $D$ such that the in-degree and out-degree of every vertex is at most $k+1$. From the directed graph $D$ we construct a bipartite graph $H=(V_1\cup V_2, E_H)$ such that $V_1$, $V_2$ are copies of $V(G^*)$, say $V_1=\{v':v\in V(G^*)\}$, $V_2=\{v'':v\in V(G^*)\}$, and we include an edge $v'w''$ in $E_H$ if and only if $(v, w)$ is an arc in $D$. The bipartite graph constructed this way has maximum degree at most $k+1$, which is equal to the size of lists on every edge (inherited from $G^*$). Hence, by Galvin's theorem, we can find a proper edge colouring of $H$ such that every edge receives a colour from its list. 
		Let us colour every edge $uv$ of $G^*$, represented by $(u,v)$ in its orientation ($D$), by the colour of $u'v''$ in $H$. Since the edge colouring of $H$ is proper and every vertex $v$ of $G^*$ is represented by its two copies in $H$, then $v$ can be incident with at most $2$ edges coloured the same in $G^*$.
The same holds for the original copy of $\overline{G}$, which is a subgraph of $G^*$.
Thus, exploiting the bijection between the edge set of $G$ and $\overline{G}$, we obtain a colouring $\omega$ of $G$ from the given lists such that every vertex $v\in V$ of degree $d(v)=2(k+1)(s-1)+i$ with $i\in[2k+2]$, i.e. represented by $s$ copies in $\overline{G}$, can be incident with at most $2s = 2\lceil d(v)/(2k+2)\rceil$ edges coloured the same.
In fact, in the case when $d(v)\equiv 1 \pmod {2k+2}$, then the first copy of $v$ in $\overline{G}$ has degree $1$, and thus is incident with at most one edge in any given colour, so overall such $v$ can be incident with 
at most $2s-1 = 2\lceil d(v)/(2k+2)\rceil-1$ edges coloured the same.
We claim that the obtained $\omega$ is a $1/k$-majority edge colouring of $G$. 
Let $v$ be any vertex of $G$. 
We need to prove that at most $\lfloor\frac{d(v)}{k}\rfloor$ edges incident with $v$ can be monochromatic. 
We consider two cases.
	\begin{case}\label{Case1}
			$\left\lfloor \frac{d(v)}{k} \right\rfloor = 2t$. Let $d(v) = 2kt+i$, $i\in\{0, \dotsc, k-1\}$. By our construction, it is enough if $2\lceil d(v)/(2k+2)\rceil\leq 2t$, which holds for  
			$d(v) \leq 2t(k+1)$. This is equivalent to $i\leq 2t$, which is true if $t\geq \frac{k-1}{2}$. Since $\delta(G)\geq 2k^2-2k$, we have $t\geq k-1$, and hence, the inequality holds.
		\end{case}
		\begin{case}
			$\left\lfloor \frac{d(v)}{k} \right\rfloor=2t+1$. Let $d(v) = 2kt+k+i$, $i\in\{0, \dotsc, k-1\}$. 
			 By our construction, it is enough if $2\lceil d(v)/(2k+2)\rceil\leq 2t$ or  $2\lceil d(v)/(2k+2)\rceil\leq 2(t+1)$ for $d(v)\equiv 1 \pmod {2k+2}$. This holds for $d(v) \leq 2t(k+1)+1$, which 
is equivalent to $k+i\leq 2t+1$. This is in turn fulfilled if $t\geq k-1$. Since $\delta(G)\geq 2k^2-2k$, the smallest degree of $v$ such that $\left\lfloor \frac{d(v)}{k} \right\rfloor$ is odd is $d(v)=2k^2-k$. Hence, $t\geq k-1$ and the desired inequality holds.
		\end{case}	
		
	\end{proof}

	As mentioned, the special case of $k=2$ within Theorem~\ref{thm_galvin}  confirms Conjecture~\ref{list_maj_conj_k2} for finite graphs.

	%{\color{red} Wydaje mi się że działa też dla nieskończonych (z tw. Tichonowa), ale pewności nie mam.
	%JP: Trzeba to ustalic i skomentowac tu i owdzie, gdzie mamy na mysli skonczone grafy, a gdzie nie...zwlaszcza, ze mozna tz przeniesc wynik dla dowolnego $k$; czy ma to sens? Ewentualnie do concluding remarks...(choc moze lepiej tu)}
	%
%	\begin{corollary}
%		Every finite graph with minimum degree at least 4 admits a majority edge colouring from lists of size 3.
%	\end{corollary}
%	
%	{\color{red} Czy umieszczac powyzsze Corollary?}

	\section{Arbitrary uniform tolerance of all colours}

	Suppose now that instead of admitting $1/k$'th fraction of incident edges in one colour, every vertex $v$ accepts at most $\alpha d(v)$ such monochromatic edges, where $\alpha$ may be any fixed real number in $(0,1)$; we shall call it \emph{tolerance}.
	 Formally, for every graph $G=(V,E)$,
	 we say that  $\omega:E\to C$ is an \emph{$\alpha$-majority edge colouring} if $d_{E_c}(v)\leq \alpha d(v)$ for every $v\in V$ and $c\in C$. 
	Such generalisation of majority edge colourings was proposed already in~\cite{majority23}, but it was not considered further except for $\alpha$ of the form $1/k$ for some integer $k\geq 2$ (i.e. $1/k$-majority edge colourings). We note here that 	the same approach as the one used to prove 
	Theorem~\ref{thm_galvin} yields the following result.

		\begin{theorem}\label{thm_galvin_alpha}
		For any integer $\ell\geq 2$ and $\alpha \in (0,1)$ such that $\alpha \ell > 1$, if a graph $G$ has a minimum degree $\delta \geq \frac{2\ell-2}{\alpha \ell-1}$,  then $G$ has an  $\alpha$-majority edge colouring from lists of size $\ell$.
		\begin{proof}
			We use a similar idea as in the proof of Theorem~\ref{thm_galvin}. Let $G$ be a graph with minimum degree $\delta \geq \frac{2\ell-2}{\alpha \ell-1}$, where $\ell$ is an integer, $\alpha\in(0,1)$ and $\alpha\ell >1$. 	For every vertex $v$ of $G$ we can write its degree as $d(v) = 2\ell(s-1)+i$, where $i\in [2\ell]$. Let us partition the neighbourhood of $v$ into $s$ disjoint sets $N_1, \dotsc, N_s$ such that $|N_1|=i$ and $|N_j|=2\ell$ for $j\geq 2$. We then construct the graph $G^*$ and the bipartite graph $H$ as described in the proof of Theorem~\ref{thm_galvin}. As in that proof the proper colouring of $H$ yields a colouring $\omega$ of $G$ from the given lists such that for every vertex $v$, the number of monochromatic edges incident with $v$ cannot exceed $2\lceil d(v)/(2\ell)\rceil$ or $2\lceil d(v)/(2\ell)\rceil-1$ if $d(v)\equiv 1 \pmod{2\ell}$.
This value is thus not greater than
$$\frac{d(v)+2\ell-2}{\ell} = \alpha d(v) \frac{1+\frac{2\ell-2}{d(v)}}{\alpha \ell} 
\leq \alpha d(v) \frac{1+\frac{2\ell-2}{\delta}}{\alpha \ell} 
\leq \alpha d(v),$$			
which means that $\omega$ is an $\alpha$-majority edge colouring of $G$.	
		 
		\end{proof}
	\end{theorem}

Note that e.g. for $\ell=k+1$ and $\alpha=1/k$, we obtain $\delta\geq 2k^2$ in Theorem~\ref{thm_galvin_alpha}, which is almost as good as the bound for $\delta$ in Theorem~\ref{thm_galvin}.

	\section{Potential further extensions of the list setting\label{Section4}}

	Let us further consider a possibly most general setting in a list variant of the edge majority problem.
		Let $G=(V,E)$ be a graph, $C$ be a set of colours and suppose every edge $e$ of $G$ is endowed with a list of colours, i.e., let $L:E\rightarrow 2^C$ be given. Assume further that 
		we are given a \emph{majority tolerance function}, 
		i.e., a function $\tau:V\times C \to (0,1)$.
		A colouring $\omega:E\to C$  is said to be a \emph{$\tau$-majority $L$-colouring}
		of $G$ if $\omega(e)\in L(e)$ for every $e\in E$ and for every vertex $v\in V$ and each colour $c\in C$,
		\begin{equation}\label{tolerance_rule}
		d_{E_c}(v)\le \tau(v,c)d(v),
		\end{equation}
		i.e. no colour $c$ can be assigned to a larger fraction of the edges incident with $v$ than the majority tolerance function admits (for the given $v$ and $c$).
Note we do not admit $\tau(v,c)=0$, as the corresponding colour $c$ would be useless, nor $\tau(v,c)=1$, as then $c$ would be always admissible (at least from the point of view of one end of an edge in whose list $c$ would appear, but by the end of this section, we shall argue that in fact we should assume symmetrical points of view for both of the ends of an edge, or more generally -- for all vertices).

In order to have any chance for existence of a $\tau$-majority $L$-colouring, we must impose some minimal assumptions concerning the list assignment. In particular, we must assume that for every edge $e$ and $v\in e$, $\sum_{c\in L(e)}\tau(v,c)\geq 1$. Otherwise we could assign to all edges incident with $v$ the same fixed list of colours which does not obey this assumption, whence preventing the existence of a desirable $\tau$-majority $L$-colouring. In fact, since $d_{E_c}(v)$ must be an \emph{integer} for every colour $c$, which we do not require from $\tau(v,c)d(v)$, we should assume that 
\begin{equation}\label{Epsilon_rule}
\sum_{c\in L(e)}\tau(v,c)\geq 1+\varepsilon
\end{equation}
for some $\varepsilon >0$.
(Even if we assumed that $\tau(v,c)d(v)$ must be an integer, then e.g. if $G$ is a $d$-regular graph with $\chi'(G)=d+1$ and $L(e)=C=[d]$ for all $e\in E$, while $\tau\equiv 1/d$, i.e. $\sum_{c\in L(e)}\tau(v,c)= 1$, then a $\tau$-majority $L$-colouring would not exist.)
  It however occurs that such assumptions are not sufficient, and do not admit proving 
  a result resembling the ones in Theorems~\ref{thm_galvin} or~\ref{thm_galvin_alpha}, 
  yielding a lower bound for $\delta(G)$ which implies that $G$ is appropriately colourable from any suitable lists,
  even if we allow $\varepsilon$ to be an arbitrarily large constant or restrict our attention to bipartite graphs exclusively.
  \begin{observation}
  For all fixed $\varepsilon,\delta > 0$, there exist a bipartite graph $G=(V,E)$ with minimum degree at least $\delta$, a list assignment $L:E\rightarrow 2^C$ and a majority tolerance function $\tau:V\times C \to (0,1)$ satisfying~(\ref{Epsilon_rule}) for every edge $e\in E$ and $v\in e$
  for which there is no
      $\tau$-majority $L$-colouring of $G$. 
  
  \begin{proof}
  Fix $\varepsilon,\delta>0$. Choose any positive integer $r \geq 1+\varepsilon$ and a small constant $\beta>0$ such that 
  \begin{equation}\label{RBeta}
  r\beta<0.5. 
  \end{equation} 
  Let $G=(V,E)$ be a complete bipartite graph $K_{n,n}$ with $n\geq \delta$, where $A,B\subseteq V$ are the two maximal independent sets of $G$ (sets of bipartition). 
  Set $L(e)=[2r]$ for every edge $e\in E$ and let for every $(v,c)\in V\times C$:
  $$\tau(v,c):=\left\{\begin{array}{rcl}
  \beta &~{\rm if}~& (v,c)\in (A\times [r]) ~~{\rm or}~~ (v,c)\in(B\times [r+1,2r])\\
  1-\beta &~{\rm if}~& (v,c)\in (A\times [r+1,2r]) ~~{\rm or}~~ (v,c)\in(B\times [r])
  \end{array}\right..
   $$
Note~(\ref{Epsilon_rule}) is fulfilled for every edge $e=uv\in E$.
Suppose $\omega$ is a $\tau$-majority $L$-colouring of $G$. 
Denote by $F$ all edges $e\in E$ with $\omega(e)\in [r]$ 
 and set $\overline{F}=E\smallsetminus F$.
Then, by~(\ref{tolerance_rule}),~(\ref{RBeta}) and the definition of $\tau$,
\begin{eqnarray}
0.5n^2&>& nr\beta n 
= \sum_{v\in A} \sum_{c\in [r]}\beta d(v) 
\geq \sum_{v\in A} \sum_{c\in [r]}d_{E_{c}}(v)
= \sum_{v\in A}d_F(v)
= |F| \nonumber\\
 &=& |E|-|\overline{F}| 
 = |E| - \sum_{v\in B}d_{\overline{F}}(v) 
 = |E| - \sum_{v\in B} \sum_{c\in [r+1,2r]}d_{E_{c}}(v) \nonumber\\
 &\geq& |E| - \sum_{v\in B} \sum_{c\in [r+1,2r]}\beta d(v) 
 = n^2-nr\beta n 
 = (1-r\beta) n^2
 > 0.5n^2, \nonumber
\end{eqnarray}
a contradiction.
  \end{proof}
\end{observation}
This observation, and many other examples, suggest
we must adopt an assumption that for any fixed colour $c\in C$,
 different vertices have the same tolerance for $c$,
i.e. a tolerance should be a function related with colours exclusively (the same for every vertex).

Finally, one should introduce one more constraint to tie up a sensible model.
In particular, it seems reasonable to assume some lower bound on the value of the tolerance function.
Otherwise, for a given graph $G$ with any maximum degree $\Delta$ one could set a constant tolerance function equal e.g. $0.5/\Delta$ and assign to every edge a list with sufficiently many colours so that~(\ref{Epsilon_rule}) holds.
However, as $(0.5/\Delta)d(v)<1$ for every vertex $v$ of $G$, assigning any colour to arbitrary edge of $G$ would result in breaching~(\ref{tolerance_rule}). In the following section we investigate such a potentially most general rational model. 
We also discuss there that the last assumption above may potentially be replaced by a slightly less intuitive constraint on a limited number of elements in every list.

	\section{General setting with diversified tolerances}

	We first specify and formalize the analysed general model.
	To that end we introduce a new notation, in order to avoid potential confusion.

Let $G=(V,E)$ be a graph and $C$ be a set of colours.
Consider a list assignment $L:E\rightarrow 2^C$ 
and so-called \emph{tolerance function} of $L$, i.e. any $\alpha:C\rightarrow (0,1)$. 
	Denote by $\min(\alpha):=\min_{c'\in C}\alpha(c')$ the \emph{minimal tolerance} of $\alpha$.
If  for every edge $e$ of $G$:
	\begin{equation}\label{1+epsilonNEW}
	 \sum\limits_{c\in L(e)} \alpha(c) \ge 1+\varepsilon 
	 \end{equation}
	for some $\varepsilon>0$, then $L$ is called an \emph{$\varepsilon$-excessive} list assignment of $G$.
	A colouring $\omega:E\to C$ is said to be an \emph{$\alpha$-majority $L$-colouring}
		of $G$ if $\omega(e)\in L(e)$ for every $e\in E$ and for each vertex $v\in V$ and any colour $c\in C$,
		\begin{equation}\label{New_tolerance_rule}
		d_{E_{c}}(v)\le \alpha(c)d(v).
		\end{equation}
	We shall prove that for any fixed $\varepsilon>0$ and $a\in(0,1)$ there is $\delta_0$ such that 
	every graph $G$ with minimum degree $\delta\ge \delta_0$ admits an $\alpha$-majority $L$-colouring of $G$ from any $\varepsilon$-excessive list assignment $L$ with a tolerance function $\alpha$ 
	satisfying 	$\min(\alpha)\geq a$,
cf. Theorem~\ref{Th:MainGeneral}. Its proof shall be probabilistic. In order to simplify our calculations and facilitate usage of a symmetric version of the Lov\'asz Local Lemma, rather than its more unwieldy general form, 
we shall again use a variant of a preliminary splitting vertices argument. 
For clarity of presentation of the main proof, we single out this initial preparatory step in the form of the following observation.
	\begin{observation}\label{Obs:splitting}
	For every fixed $\varepsilon>0$, $a\in(0,1)$ and an integer $\delta_0>0$, in order to prove that
	for every graph $G$ 
	with minimum degree $\delta(G)\ge \delta_0$, any $\varepsilon$-excessive list assignment $L$ of $G$ 
	and every tolerance function $\alpha$ of $L$ 
	with $\min(\alpha)\geq a$, there is an $\alpha$-majority $L$-colouring of $G$,
	it is enough to prove this statement for graphs $G$ with all vertex degrees in $[\delta_0,2\delta_0-1]$.
	\begin{proof}
			Let $G$ be a graph with minimum degree $\delta(G)$ not smaller than $\delta_0$ (which satisfies all graphs with degrees in $[\delta_0,2\delta_0-1]$) and let $L$ be any $\varepsilon$-excessive list assignment of $G$ with tolerance function $\alpha$ such that $\min(\alpha)\geq a$. Let $v$ be a vertex of $G$, i.e. $d_G(v) \geq \delta_0$. 	Similarly as in the proofs of Theorems~\ref{thm_galvin} and~\ref{thm_galvin_alpha}, we may split $v$ by first partitioning the neighbourhood of $v$ into $N_1,\ldots,N_s$ with $\delta_0\leq |N_i|\leq 2\delta_0-1$, and then associating each $N_i$ as a neighbourhood to a new copy of $v$, say $v_i$ (and deleting $v$ itself).
This operation yields a bijection between the edges of $G$ and the edges of the resulting graph.
After processing and possibly splitting all vertices of $G$, we end up with a graph $\overline{G}=(\overline{V},\overline{E})$ with  all vertex degrees in $[\delta_0,2\delta_0-1]$, naturally inheriting the $\varepsilon$-excessive list assignment from $G$, due to the mentioned bijection between the edges of $G$ and $\overline{G}$.
By our assumption, there is an $\alpha$-majority colouring of $\overline{G}$ from these lists,
i.e. such that for every given colour $c$ and a copy $v_i$ in $\overline{G}$ of any vertex $v$ of $G$, we have $d_{\overline{E}_c}(v_i)\leq \alpha(c) d_{\overline{G}}(v_i)$. Therefore, if $v$ was split into $s$ copies, then the bijectively mapped edge colouring from $\overline{E}$ to $E$ fulfils: 
$$d_{E_c}(v) = \sum_{i=1}^s d_{\overline{E}_c}(v_i) \leq \sum_{i=1}^s \alpha(c)d_{\overline{G}}(v_i) = \alpha(c)\sum_{i=1}^s d_{\overline{G}}(v_i) = \alpha(c) d_G(v).$$
Therefore,  the resulting $L$-colouring of $G$ is an $\alpha$-majority colouring.

		\end{proof}
	\end{observation}

	Let us also recall two standard tools applicable within the probabilistic method, see e.g.~\cite{LLL} and~\cite{RandOm}, respectively.
	
	\begin{lemma}[Lov\'asz Local Lemma]\label{LLL}
Let $\Omega$ be a finite family of events in any probability space. Suppose that every event $A \in \Omega $ is mutually independent of a set of all the other events in $\Omega$ but at most $D$, and that $\mathbf{Pr}(A) \leq p$ for each $A \in \Omega$. If 
$$ep(D + 1) \leq 1,$$
then $\mathbf{Pr}\left(\bigcap_{A \in \Omega} \overline{A}\right) > 0$. 
\end{lemma}

\begin{lemma}[Chernoff Bound]\label{Ch} Let $X= \sum_{i=1}^n X_i$, where $X_i=1$ with probability $p_i$, $X_i = 0$ with probability $1 - p_i$ and all $X_i$ are independent. Then, for every $0<t\leq \mathbf{E}(X)$,
$$ \mathbf{Pr}\left(X > \mathbf{E}(X)+t\right)< \exp\left(-\frac{t^2}{3\mathbf{E}(X)}\right).$$
\end{lemma}
One may easily verify that the Chernoff Bound above is also applicable in the case when only an upper bound for $\mathbf{E}(X)$ is known, say $\mathbf{E}(X)\leq \sigma$, whence Lemma~\ref{Ch} implies that $ \mathbf{Pr}(X > \sigma+t)< \exp(-t^2/(3\sigma))$ for every $0<t\leq\sigma$.
	
	We shall now prove the theorem aforementioned at the beginning of this section. In order to provide a relatively uncomplicated formula for $\delta_0$, we impose some mild assumptions on $\varepsilon$, that is we prove the theorem for $\varepsilon$ only slightly below $1$, i.e. $\varepsilon\leq 0.9$. Note at the same time that $\delta_0$ derived from Theorem~\ref{Th:MainGeneral} for $\varepsilon=0.9$ shall be the more sufficient for larger values of $\varepsilon$, by formulation of our problem itself.
	\begin{theorem}\label{Th:MainGeneral}
	For every fixed $\varepsilon\in(0,0.9]$, $a\in(0,1)$, for every graph $G$ 
	with minimum degree $\delta\ge \lceil 626a^{-1}\varepsilon^{-2}\ln(a\varepsilon)^{-1}\rceil$ 
	and each $\varepsilon$-excessive list assignment $L$ of $G$ 
	associated with any tolerance function $\alpha$ of $L$ 
	with $\min(\alpha)\geq a$, there is an $\alpha$-majority $L$-colouring of $G$.
	
	\begin{proof}
	Let us fix $\varepsilon\in(0,0.9]$, $a\in(0,1)$. Set
	\begin{equation}\label{delta_0_def}
	\delta_0 = \left\lceil 626a^{-1}\varepsilon^{-2}\ln(a\varepsilon)^{-1}\right\rceil.
	\end{equation}
	By Observation~\ref{Obs:splitting} we may assume that $G=(V,E)$ is a graph with minimum degree 
	$\delta\geq \delta_0$ and maximum degree  $\Delta\leq 2\delta_0-1$.	
	Fix any $\varepsilon$-excessive list assignment $L:E\rightarrow 2^C$ of $G$  
	and a tolerance function $\alpha:C\rightarrow (0,1)$  of $L$
	with $\min(\alpha)\geq a$.
	
	From each list $L(e)$ we first remove some colours, if necessary, so that afterwards~(\ref{1+epsilonNEW}) still holds and:
	\begin{equation}\label{Bounded_list_size}
	|L(e)|\leq \left\lceil\frac{1+\varepsilon}{a}\right\rceil < \frac{1+\varepsilon+a}{a}. 
	\end{equation}
	For each resulting list $L(e)$ independently we shall randomly choose one of its colour, with probability for each colour $c$ proportional to the value of the tolerance function: $\alpha(c)$, and show that with positive probability the resulting edge colouring meets our requirements. 
	Formally, for every edge $e\in E$ we define below a random variable $X_e$ accountable for prescribing a colour to $e$. For the sake of symmetrization of notation, we extend its natural domain (sample space) of $L(e)$ to entire $C$ (prescribing probability $0$ to colours outside $L(e)$). Thus, for every $e\in E$ and $c\in C$, the variable $X_e$ takes value $c$ with probability:
	\begin{equation}\label{pec}
	p_{e,c}:= \left\{\begin{array}{cll}\frac{\alpha(c)}{\sum_{c'\in L(e)}\alpha(c')}&~{\rm if}~& c\in L(e)\\
	0&~{\rm if}~& c\notin L(e)
	\end{array}\right..
	\end{equation}
	For each $e\in E$ and $c\in C$, we further define the binary random variable $Y_{e,c}$ taking value $1$ if $X_e=c$, and $0$ otherwise. Note that since $L$ is $\varepsilon$-excessive, i.e.~(\ref{1+epsilonNEW}) holds, then by~(\ref{pec}), we obtain that:
	\begin{equation}\label{PrYec1}
	\mathbf{Pr}\left(Y_{e,c}=1\right)=p_{e,c}\leq \frac{\alpha(c)}{1+\varepsilon}.
	\end{equation}
	Denote by $L_v$ the set of all colours appearing in lists incident with a vertex $v$, i.e., $L_v=\bigcup_{e\in E(v)}L(e)$. For every $c\in L_v$ we finally define the following random variable:
	$$Z_{v,c}=\sum_{e\in E(v)}Y_{e,c}.$$
	Note that all random variables $Y_{e,c}$ in the sum above are independent, thus the Chernoff Bound applies to $Z_{v,c}$. Note also that $Z_{v,c}$ represents the number of edges incident with $v$ which are coloured $c$ within our random edge colouring. In order to prove that the colouring meets our expectations 
	 it is thus sufficient to prove that (with positive probability) for no $v\in V$ and $c\in L_v$ the following event occurs: 
	$$A_{v,c}:~Z_{v,c}>\alpha(c)d(v).$$
	By~(\ref{PrYec1}),
	$$\mathbf{E}\left(Z_{v,c}\right)=\sum_{e\in E(v)}\mathbf{E}\left(Y_{e,c}\right)
	=\sum_{e\in E(v)}p_{e,c}
	\leq \frac{\alpha(c)d(v)}{1+\varepsilon}.$$
	Thus, by the Chernoff Bound, 
	\begin{eqnarray}
	\mathbf{Pr}\left(A_{v,c}\right)
	&=& \mathbf{Pr}\left(Z_{v,c}>\frac{\alpha(c)d(v)}{1+\varepsilon}+\frac{\varepsilon\alpha(c)d(v)}{1+\varepsilon}\right) \nonumber\\
	&\leq& \exp\left(-\frac{\varepsilon^2\alpha(c)d(v)}{3(1+\varepsilon)}\right) \nonumber\\
	 &\leq& \exp\left(-\frac{\varepsilon^2ad(v)}{3(1+\varepsilon)}\right).  \label{PrAvcBound}
	 \end{eqnarray}
	 Set 
	 $$b=\frac{\varepsilon^2a}{3(1+\varepsilon)},\hspace{1cm} 
	 \gamma=626a^{-1}\varepsilon^{-2}\ln(a\varepsilon)^{-1},$$ 
	 hence 
	 \begin{equation}\label{GammaSmall}
	 \gamma\leq \lceil\gamma\rceil = \delta_0 \leq \delta \leq d(v)
	 \end{equation}
	 for every vertex $v\in V$.
	 Consider the function 
	 \begin{equation}\label{Fdef}
	 f(x)=\frac{2e(1+\varepsilon+a)}{a}x^2e^{-bx}
	 \end{equation}
	 and note that
	 $f'(x)=\frac{2e(1+\varepsilon+a)}{a}xe^{-bx}(2-bx)$, and hence $f$ is decreasing for $x\geq \frac{2}{b}$, while 
	 $$\gamma=\frac{626\ln(a\varepsilon)^{-1}}{3(1+\varepsilon)b}>\frac{626\ln(0.9)^{-1}}{5.7b}>\frac{2}{b}.$$
	 Therefore, $f$ is decreasing for $x\geq \gamma$. Hence, for every such $x$: 
\begin{eqnarray}
	 f(x)\leq f(\gamma) &=& \frac{2e(1+\varepsilon+a)}{a}626^2a^{-2}\varepsilon^{-4}\ln^2(a\varepsilon)^{-1}e^{-\frac{626\ln(a\varepsilon)^{-1}}{3(1+\varepsilon)}} \nonumber\\
	 &\leq&  2e(1+\varepsilon+a)626^2(a\varepsilon)^{-4}\ln^2(a\varepsilon)^{-1}(a\varepsilon)^{\frac{626}{3(1+\varepsilon)}} \nonumber\\
%	 	 &=&  2e(1+\varepsilon+a)626^2\ln^2(a\varepsilon)^{-1}(a\varepsilon)^{\frac{626}{3(1+\varepsilon)}-4} \nonumber\\
	 	 &\leq&  5.8e626^2\ln^2(a\varepsilon)^{-1}(a\varepsilon)^{\frac{626}{5.7}-4} \nonumber\\
	 &\leq& 5.8e626^2 \ln^2(0.9)^{-1}0.9^{\frac{626}{5.7}-4}  < 1, \label{fxIneq}
\end{eqnarray}
where the last weak inequality above follows by the fact that $f_1(x)=\ln^2x^{-1}\cdot x^{\frac{626}{5.7}-4}$ is increasing for $x$ between $0$ and $\exp(-\frac{2}{\frac{626}{5.7}-4})>0.9$. 
Let us aggregate all $A_{v,c}$ related with a given vertex $v$ within a single undesirable event:
$$B_v:~~{\rm therere~exist~} c\in L_v {\rm ~for~which~} Z_{v,c}>\alpha(c)d(v).$$
By~(\ref{PrAvcBound}), (\ref{Bounded_list_size}), (\ref{Fdef}), (\ref{fxIneq}) and~(\ref{GammaSmall}),
\begin{eqnarray}
\mathbf{Pr}(B_v) &\leq& |L_v|\exp\left(-\frac{\varepsilon^2ad(v)}{3(1+\varepsilon)}\right) 
< d(v)\frac{1+\varepsilon+a}{a} e^{-b\cdot d(v)} \nonumber\\
&=&\frac{f(d(v))}{2e\cdot d(v)} 
< \frac{1}{2e\cdot d(v)}  
\leq \frac{1}{2e\cdot \delta_0}.  \label{PrBvEstimate}
\end{eqnarray}
Denote $p=(2e\cdot\delta_0)^{-1}$, hence by~(\ref{PrBvEstimate}), $\mathbf{Pr}(B_v)<p$ for every $v\in V$. Since every event $B_v$ is mutually independent of all but at most $D=d(v)\leq 2\delta_0-1$ other events $B_{v'}$, by the Lov\'asz Local Lemma, there is a choice of colours from the lists so that none of the events $B_v$, $v\in V$, holds.	
	\end{proof}
	\end{theorem}

	As mentioned earlier, instead of imposing a lower bound for the tolerance function $\alpha$, one may analyse a somewhat related condition where we admit a limited number of elements in each list. In such a case, the following result may be derived from Theorem \ref{thm_galvin}.
	
	Let $\ell\geq 2$ be an integer and let $\varepsilon>0$. We say that a graph $G$ is $(\ell, \varepsilon)$-majority edge-choosable if for each tolerance function $\alpha:C\to (0,1)$ and every collection of lists such that  $|L(e)|\leq \ell$ and $\sum_{c\in L(e)} \alpha(c) \geq 1 + \varepsilon$ for every edge $e$ of $G$, there is an $\alpha$-majority $L$-colouring of $G$. 
	
	\begin{theorem}\label{tw14}
			For any integer $\ell\geq 2$ and any real number $\varepsilon>0$, there exists an integer $\delta_{(\ell, \varepsilon)}$ such that for every graph $G$ with minimum degree $\delta \geq \delta_{(\ell, \varepsilon)}$, each tolerance function $\alpha:C\to (0,1)$ and every
			$\varepsilon$-excessive list assignment $L$ of $G$
			such that  $|L(e)|\leq \ell$ for every edge $e$ of $G$, there is an $\alpha$-majority $L$-colouring of $G$. It is sufficient to let $\delta_{(\ell, \varepsilon)} = 2\left\lceil \frac{\ell}{\varepsilon} \right\rceil^2-2\left\lceil \frac{\ell}{\varepsilon} \right\rceil$.
		\begin{proof}
	Consider a graph $G=(V,E)$. Fix $\alpha:C\to(0,1)$	and $L:E\to 2^C$ such that 
	$|L(e)|\leq \ell$ and $\sum_{c\in L(e)} \alpha(c) \geq 1 + \varepsilon$ 
			for every edge $e\in E$.
			Set $k:=\lceil \ell/\varepsilon \rceil$. To every edge $e$ of $G$ we assign a modified  list of colours $L'(e)$ constructed of $L(e)$ as follows.
For every colour $c \in C$ we replace $c$ with $s=\lfloor k\alpha(c) \rfloor$ copies of $c$, say $c_1,\ldots,c_s$, and we set $\alpha'(c_i)=1/k$ for each $i\in[s]$. Then, 
\begin{equation}\label{C_i-C-frequencies}
\sum_{i=1}^s\alpha'(c_i)=\frac{s}{k} > \alpha(c)-\frac{1}{k}.
\end{equation}
			 Note also that if $G$ has an $\alpha'$-majority $L'$-colouring $\omega'$,
			 then it yields an $\alpha$-majority $L$-colouring $\omega$ of $G$ where we associate a colour $c$ to every edge $e$ with any copy $c_i$ of $c$ assigned by $\omega'$. Indeed, it is enough to observe that for every vertex $v$ and a colour $c$ with $s=\lfloor k\alpha(c) \rfloor$, the number of edges in $E(v)$ coloured $c$ by $\omega$ can be expressed by the following sum related with frequencies of the copies of $c$ under $\omega'$: 
			$$d_c(v)=\sum_{i=1}^sd_{c_i}(v)\leq \sum_{i=1}^s\alpha'(c_i)%\frac{1}{k}
			d(v) = s\frac{1}{k}d(v)
			\leq \alpha(c)d(v).$$
	On the other hand, by~(\ref{C_i-C-frequencies}), for each $e\in E$,	
			\[ \sum_{c\in L(e)} \alpha(c) - \sum_{c'\in L'(e)} \alpha'(c') < \frac{\ell}{k} \leq \varepsilon, \]
			hence 
			$\sum_{c'\in L'(e)} \alpha'(c') > 1$. 
			Consequently, as $\alpha'(c')=\frac{1}{k}$ for $c'\in L'(e)$,
			every $L'(e)$ must be of size at least $k+1$. Thus, by Theorem \ref{thm_galvin} we obtain that if $\delta(G) \geq 2k^2-2k$, then $G$ has a $1/k$-majority edge colouring from the lists $L'(e)$, $e\in E$, which simply is an $\alpha'$-majority $L'$-colouring, guaranteeing the existence of a desired $\omega$.
		\end{proof}
	\end{theorem}

Now, going back to our initial setting, with $\min(\alpha)\geq a$,  similarly as within the proof of Theorem~\ref{Th:MainGeneral}, 
one may always remove some excessive colours from each $L(e)$ to assure that $|L(e)|\leq\lceil\frac{1+\varepsilon}{a}\rceil$ for every edge $e$. Thus, Theorem~\ref{tw14} implies an upper bound (for the optimal value of a threshold value of $\delta$) 
of the form $O((a\varepsilon)^{-2})$ (for a bounded $\varepsilon$, e.g., $\varepsilon\leq 0.9$), which is in general worse than the 
bound of order $O(a^{-1}\varepsilon^{-2}\ln(a\varepsilon)^{-1})$, stemming from Theorem~\ref{Th:MainGeneral},
with respect to $a$, though slightly better in terms of $\varepsilon$.
However, relating such a bound with $\ell$, rather than $a$, might be frequently quite beneficial, e.g. when only a few colours in the lists have very small tolerance, while the majority have large (whence $\ell$ may be very small compared to $a^{-1}$). 

In the following section we analyse one more special model of lists assignment, where we impose the same vector of tolerances of colours in each list. It occurs that in such a setting, one may
provide a bound of order $O(\ell\varepsilon^{-2}\ln(\ell\varepsilon)^{-1})$, which is not a function of the  less preferred parameter $a$, and is better
from the discussed bound from Theorem~\ref{tw14} with respect to $\ell$ (at a miniscule cost in view of $\varepsilon$). 
At the same time, it provides a better multiplicative constant than the one resulting from Theorem~\ref{Th:MainGeneral},
and is applicable within a general non-list setting, cf. 
Section~\ref{Concluding_section}.

	\section{Uniform edge vector of tolerances}

Suppose we want to symmetrize our setting slightly by imposing a relatively natural restriction (usually adapted for regular list colourings, where the lists are typically uniform with respect to their structure, in particular sizes).
Namely, we still admit different lists of colours associated to distinct edges, whereas every colour may also have a different tolerance associated.
	However, we shall require that each list consists in the same number of colours of every given \emph{type}, where two colours are considered the same type if they have the same tolerance. 
	Note this holds in particular in a general non-list setting of majority colourings, with different tolerances permissible for distinct colours; more comments on this issue are included in the last section.

	Let $G=(V,E)$ be a graph and let $\Lambda = (\alpha_1, \alpha_2, \ldots, \alpha_\ell)$ be a vector of real numbers (not necessarily pairwise distinct) such that $\alpha_i \in (0,1)$ for each $i$ and 
	\begin{equation}\label{AlphaiEps}
	\sum_{i\in [\ell]}\alpha_i = 1+\varepsilon
	\end{equation}
	for some $\varepsilon>0$. We call such $\Lambda$ an \emph{$\varepsilon$-excessive} \emph{tolerance vector}.  Let $C$ be a set of colours and $\alpha:C\rightarrow (0,1)$ be a function assigning to every colour from $C$ some weight $\alpha_i$ from the vector $\Lambda$; we call it a \emph{$\Lambda$-tolerance function}. For every edge $e$ of $G$ let $L(e) = (c_1, \dotsc, c_\ell)$ be a list of colours such that $c_i \in C$ and $\alpha(c_i)=\alpha_i$ for $i\in[\ell]$. We call such list a 
	\emph{$\Lambda$-list}. A \emph{$\Lambda$-list assignment} of $G$ is an $L:E\to 2^C$ such that $L(e)$ is a $\Lambda$-list for every $e\in E$. 
	
	To investigate this setting we shall apply a probabilistic approach similar to the one used within the proof of Theorem~\ref{Th:MainGeneral}, that is, we shall randomly choose colours from the lists with probabilities (nonlinearly) proportional to their tolerances. We shall be able to introduce at least two important optimising refinements to our reasoning, though. One of these relies on a technical lemma concerning 
	conditions implying maximisation of probability of a certain type of events, see Lemma~\ref{The_same_lists}  below.
	
	Before we present this lemma, we first formulate a rather obvious auxiliary observation, which states that the sum of
	a fixed number of consecutive elements from any specific row of Pascal's triangle is the larger the more central the set of these elements is. We also include the proof of this observation,
	 for the sake of completeness. 
For technical reasons we admit within the following observations binomial coefficients ${n \choose m}$ with $m<0$ or $m>n$, which by definition equal $0$.
	
	\begin{observation}\label{ObsPascal}
	For any positive integer $z$ and integers $a,b,c,d$ %$a,b,c,d\in[0,z]$ 
	with $b-a=d-c\geq 0$ and %$\min\{a,z-b\}\leq \min\{c,z-d\}$,  
	$|\frac{a+b}{2}-\frac{z}{2}|\geq |\frac{c+d}{2}-\frac{z}{2}|$,
	\begin{equation}\label{ChooseThesis}
	\sum_{i=a}^b{z\choose i} \leq \sum_{i=c}^d{z\choose i}. 
	\end{equation}
	
\begin{proof}
We shall use the following straightforward claim on 
integers $k,k'$.
%nonnegative integers $k,k'$ not exceeding $z$. 
\begin{claim}\label{ClaimChoose}
If ${z\choose k} \leq {z\choose k'}$, %and $k\leq \frac{z}{2}$, 
$k\leq \frac{z}{2}$ and $k\leq k'$,
then ${z\choose k-1} \leq {z\choose k'+1}$.
\end{claim}

We may assume $[a,b]$ and $[c,d]$ are disjoint, as the common summands on both sides of~(\ref{ChooseThesis}) do not influence the inequality, while removing them results in modified values of $a,b,c,d$, still fulfilling all assumptions of the observation.
We may further assume that $a\leq c$. (If this is not the case, it suffices to exchange every ${z\choose x}$ with 
${z\choose z-x}$ in~(\ref{ChooseThesis}).)
Hence, $b\leq \frac{z}{2}$ and $b<c\leq z-b$, where the weak inequalities follow by the observation's assumptions on a more central localisation of $[c,d]$ than $[a,b]$ with respect to $[0,z]$. Thus, ${z\choose b}\leq {z\choose c}$.
Consequently, by Claim~\ref{ClaimChoose}, ${z\choose b-1}\leq {z\choose c+1}$ (where $b-1\leq \frac{z}{2}$ and $b-1<c+1$). By analogous repeated applications of Claim~\ref{ClaimChoose}, we likewise obtain that ${z\choose b-j}\leq {z\choose c+j}$ for $j\in[0,b-a]$, and hence~(\ref{ChooseThesis}) follows.  
\end{proof}
	\end{observation}
	
	\begin{lemma}\label{The_same_lists}
	Let $\Lambda = (\alpha_1, \alpha_2, \ldots, \alpha_\ell)$  be an $\varepsilon$-excessive tolerance vector
	and let $L:E\to 2^C$ be a \emph{$\Lambda$-list assignment} of a graph $G=(V,E)$, whereas $\alpha:C\rightarrow (0,1)$ is the corresponding $\Lambda$-tolerance function.
	Let $v\in V$ be a fixed vertex with $E(v)=\{e_1,e_2,\ldots,e_d\}$
	and let $p_1,\ldots,p_\ell\geq 0$ be constants such that $\sum_{i\in [\ell]}p_i=1$, $p_i< \alpha_i$ for $i\in[\ell]$ and $p_i=p_j$ whenever $\alpha_i=\alpha_j$. 
	Suppose that for each $e\in E$ we independently choose 
	a colour from its list $L(e)=(c_1,c_2,\ldots,c_\ell)$ at random where each colour $c_i$ is drawn with probability $p_i$. 
	Then the probability that for some colour $c\in C$ we have $d_{E_c}(v)>\alpha(c)d$ for the fixed $v$ is the largest in the case when all lists $L(e_1),L(e_2),\ldots,L(e_d)$ are the same.
	
	\begin{proof}
	Fix any $G=(V,E)$, $\Lambda = (\alpha_1, \ldots, \alpha_\ell)$, $\alpha:C\rightarrow (0,1)$ and $p_i$, $i\in[\ell]$ consistent with the lemma's assumptions. Let $L_1,\ldots,L_d$ be any $\Lambda$-lists associated to pairwise distinct edges $e_1,\ldots,e_d$, resp.,  incident with a vertex $v$ of degree $d$ in $G$.
	We call  $\mathcal{L}=\{L_1,\ldots,L_d\}$  a \emph{$\Lambda$-family of $v$}.
	Let 
	$\bigcup \mathcal{L} = \bigcup_{i\in [d]} L_i=\{a_1,\dotsc,a_n\}\subseteq C$ denote the set of all colours represented within $\mathcal{L}$; note $n\geq \ell$. 
	Let further $\# a_i$ denote the number of occurrences of the colour $a_i$ in the lists in $\mathcal{L}$.
	Then by $W(\mathcal{L}) = (w_1,\dotsc,w_n)$ we denote the nonincreasingly sorted vector of $\# a_i$, $i\in[n]$, and call it the \emph{vector of colour occurrences} in $\mathcal{L}$. 
			
			Let us introduce an order in the set of all possible $\Lambda$-families of $v$ (including colours from $C$) consistent with the lexicographic order within the set of the corresponding vectors of colour occurrences. 
			To be precise, for any given $\Lambda$-families  $\mathcal{L}_1, \mathcal{L}_2$ of $v$ with $W(\mathcal{L}_1)=(w_1,\dotsc,w_n)$, $W(\mathcal{L}_2)=(u_1,\dotsc,u_n)$ (if these vectors are not equal length, we append zeroes to the shorter one to match their lengths), we assume 
			$\mathcal{L}_1 \prec \mathcal{L}_2$ 
			if for some $i$ we have $w_i < u_i$ 
			and $w_j = u_j$ for all $j<i$. Hence, `$\prec$' 
			defines a strict partial order, while every maximal element with respect to $\prec$ in the set of all possible $\Lambda$-families of $v$  consists in uniform lists associated to all edges in $E(v)$
			(i.e. with exactly the same $L$ associated to every edge incident with $v$, where $L$ is any $\Lambda$-list).
			
			In order to prove the lemma it is thus sufficient to show that for any $\Lambda$-family $\mathcal{L}$ of $v$ which is not a maximal element with respect to  $\prec$, there exists a $\Lambda$-family $\mathcal{L}'$ of $v$ such that $\mathcal{L} \prec \mathcal{L}'$ and the probability that for some colour $c\in C$ we have $d_{E_c}(v)>\alpha(c)d$  when colours incident with $v$ are chosen from the lists in $\mathcal{L}$ does not decrease if these colours are chosen from the lists in $\mathcal{L}'$ instead.
			
			Suppose then that $\mathcal{L} = \{L_1,\dotsc,L_d\}$ is a non-maximal $\Lambda$-family $\mathcal{L}$ of $v$.
			Set $\bigcup \mathcal{L} = \{a_1,\dotsc,a_n\}$ and let $a_s$ be the least frequently appearing colour in $\mathcal{L}$, i.e. 
			$\# a_s = \min_i (\# a_i)$. Since  $\mathcal{L}$ is not maximal,  
			we have $1 \leq \# a_s < d$. Thus, $a_s$ must appear in some list, say $L_{i_1}$
			 and not appear in some other list, say $L_{i_2}$  in $\mathcal{L}$.
			 However, as $L_{i_1},L_{i_2}$ are $\Lambda$-lists, where $\Lambda=(\alpha_1, \ldots, \alpha_\ell)$ is fixed, both of the lists must contain the same number of colours $a$ with $\alpha(a)=\alpha(a_s)$. Therefore, there must exist $a_t\in L_{i_2}$ with $\alpha(a_t)=\alpha(a_s)$ such that $a_t\notin L_{i_1}$.			  
			Let $I_{s-t}$ be the set of all indices $i\in[d]$ such that $a_s \in L_i$ and  $a_t \notin L_i$, hence in particular, $i_1\in I_{s-t}$.
			
			Let $\mathcal{L}'=\{L'_1,\ldots,L'_d\}$ be defined by setting:
			\begin{equation}\label{L'def}
			L'_i=\left\{\begin{array}{ll}
			L_i\cup \{a_t\}\smallsetminus \{a_s\}, &~~{\rm if}~~~~  i\in I_{s-t},\\
			L_i, &~~{\rm otherwise.}
			\end{array}\right.
			\end{equation}
			Note that $\mathcal{L}'$ is a $\Lambda$-family of $v$ and, by the choice of $a_s$, $\mathcal{L} \prec \mathcal{L}'$.  
			Denote the following event:
			$$A: ~{\rm there~ exists}~ c\in C ~{\rm such~ that}~ d_{E_c}(v)>\alpha(c)d.$$
			We shall show that the probability that $A$ occurs does not decrease when colours incident with $v$ are drawn from lists in $\mathcal{L}'$ instead of $\mathcal{L}$, thus finishing the proof. 
For the sake of clarity, we shall use the notation  $\mathbf{Pr_1}$ 
to denote the probability in the setting when the colours are chosen from lists in $\mathcal{L}'$, and the
notation $\mathbf{Pr_0}$ 
in the case when the colours are drawn from lists in $\mathcal{L}$. We shall thus show that
\begin{equation}\label{P0P1ineq1}
\mathbf{Pr_0}(A)\leq \mathbf{Pr_1}(A).
\end{equation}
			Note we may 
			bijectively couple each colouring of $E(v)$, say $\omega_v$ chosen from $\mathcal{L}$ with a colouring $\omega'_v$ chosen from $\mathcal{L}'$, by recolouring to $a_t$ each edge $e_i$ with $i\in I_{s-t}$ which is coloured $a_s$ under $\omega_v$. 
				Note moreover that as $\alpha(a_t)=\alpha(a_s)$, then the probability 
				of choosing $a_s$ from such a list $L_i$ is the same as the probability of choosing $a_t$ from the corresponding list $L'_i$ (while all colour choices are independent). 
			Thus, the colourings $\omega_v,\omega'_v$ shall be chosen with the same probabilities, 
			i.e. $\mathbf{Pr_0}(\omega_v)=\mathbf{Pr_1}(\omega'_v)$.
			
			Let 
			$A_J$ be an event such that for some colour $c\in C$, the number of times the colour $c$ was chosen from the given lists $L_k$ (or $L'_k$) with indices $k\in J$ is greater than $\alpha(c)d$. In particular, $A_{[d]}=A$. Let 
			$I_{s,t}$ be a set of all indices $k$ such that $a_s\in L_k$ or $a_t \in L_k$ (or equivalently: $a_s\in L'_k$ or $a_t \in L'_k$). Let $I$ denote the (random) set of indices such that $k\in I$ iff the colour chosen from the list $L_k$ 
			 (or $L'_k$, through the mentioned bijective coupling)
			  is either $a_s$ or $a_t$, thus $I\subseteq I_{s,t}$. 
			
			By the law of total probability, for $i\in\{0,1\}$: 
			\[ \mathbf{Pr_i} (A) = \sum\limits_{I'\subseteq I_{s,t}} 
			\mathbf{Pr_i}(I=I')\cdot \mathbf{Pr_i}(A | I=I'). \]
	Let us fix any $I'\subseteq I_{s,t}$. 
			Note that by the bijective coupling and comments above, $\mathbf{Pr_0}(I=I')=\mathbf{Pr_1}(I=I')$. 
			In order to show~(\ref{P0P1ineq1}) 
			it thus suffices to prove that 
			\begin{equation}\label{P0P1ineq2}
			\mathbf{Pr_0}(A | I=I')\leq \mathbf{Pr_1}(A | I=I').
			\end{equation}
			Note that for $i\in\{0,1\}$:
\begin{equation}\label{P0P1ineq3}
\mathbf{Pr_i}(A | I=I') = \mathbf{Pr_i}(A_{[d]\smallsetminus I'}|I=I') + \mathbf{Pr_i}(\overline{A_{[d]\smallsetminus I'}} \cap A|I=I'). 
\end{equation}
			Note further that no edge $e_i$ with $i\in [d]\smallsetminus I'$ can have colour $a_s$ or $a_t$ assigned, given that $I=I'$. However, apart from $a_s$ and $a_t$, the lists in $\mathcal{L}$ and $\mathcal{L}'$ coincide. Thus,
			\begin{equation}\label{P0P1ineq4}
			\mathbf{Pr_0}(A_{[d]\smallsetminus I'}|I=I')=\mathbf{Pr_1}(A_{[d]\smallsetminus I'}|I=I').
			\end{equation} 
			Furthermore, 
			given that $I=I'$, the
			set of colours chosen from the lists with indices in $[d]\smallsetminus I'$
			is disjoint from the set of colours picked from the lists with indices in $I'$.
			Therefore, for $i\in\{0,1\}$:
	%{\color{red}
			\begin{eqnarray}
			\mathbf{Pr_i}(\overline{A_{[d]\smallsetminus I'}} \cap A|I=I') 
			&=& \mathbf{Pr_i}(\overline{A_{[d]\smallsetminus I'}} \cap A_{I'}|I=I') \nonumber\\
			&=& \mathbf{Pr_i}(\overline{A_{[d]\smallsetminus I'}}|I=I') \cdot \mathbf{Pr_i}(A_{I'}|I=I'), \label{P0P1ineq5}
			\end{eqnarray} 
			where the last equality follows by the fact that choices of colours for distinct edges are independent. Since the lists in $\mathcal{L}$ and $\mathcal{L}'$ coincide apart from $a_s$ and $a_t$, we have
			\begin{equation}\label{P0P1ineq6}
				\mathbf{Pr_0}(\overline{A_{[d]\smallsetminus I'}}|I=I')=\mathbf{Pr_1}(\overline{A_{[d]\smallsetminus I'}}|I=I').
			\end{equation} 
		
			By (\ref{P0P1ineq3}), (\ref{P0P1ineq4}), (\ref{P0P1ineq5}) and (\ref{P0P1ineq6}), 
	%}
			in order to show (\ref{P0P1ineq2}) (hence also (\ref{P0P1ineq1})), it suffices to prove that 
\begin{equation}\label{P0P1ineq7}
\mathbf{Pr_0}(A_{I'}|I=I')\leq \mathbf{Pr_1}(A_{I'}|I=I'). 
 \end{equation}

			Since the probabilities of choosing $a_s$ and $a_t$ from any given list including these colours are the same,  and hence  (provided that $I=I'$), any attainable choice of $\{a_s,a_t\}$-colouring of the edges in  $E':=\{e_i:i\in I'\}$ from the lists in $\mathcal{L}$ or $\mathcal{L}'$ is equally probable, 
			 in order to prove~(\ref{P0P1ineq7}), it suffices to show that the number of 
			 $\{a_s,a_t\}$-colourings of $E'$ 
			  for which $d_{E'_{a_s}}(v)> \alpha(a_s)d$ or $d_{E'_{a_t}}(v)> \alpha(a_s)d$ 
			 from appropriate lists in  $\mathcal{L}'$ is not smaller than from the corresponding lists in  $\mathcal{L}$. Denote the set of all such colourings from  $\mathcal{L}'$ by $\Omega'$, and the set of all such colourings from $\mathcal{L}$ by $\Omega$.
			
			Set $r=\lfloor \alpha(a_s)d +1 \rfloor$, hence $r>\alpha(a_s)d$.
			Let us denote cardinalities of three subsets partitioning the set of
			all lists $L_i$ with $i\in I'$:
			\begin{align}
			x&=|\{L_i:i\in I', L_i\cap\{a_s,a_t\}=\{a_t\}\}| \label{xDef}\\
			y&=|\{L_i:i\in I', L_i\cap\{a_s,a_t\}=\{a_s\}\}| \label{yDef}\\
			z&=|\{L_i:i\in I', L_i\cap\{a_s,a_t\}=\{a_s,a_t\}\}| \label{zDef}
			\end{align}
	Thus, by~(\ref{L'def}), 
			\begin{align}
			x+y&=|\{L'_i:i\in I', L'_i\cap\{a_s,a_t\}=\{a_t\}\}| \label{x+yDef'}\\
			z&=|\{L'_i:i\in I', L'_i\cap\{a_s,a_t\}=\{a_s,a_t\}\}| \label{zDef'} 
			\end{align}
			(Note that by construction, none of the lists $L'_i$ with $i\in I'$ contains $a_s$ but does not contain $a_t$.)
			Note that if $x+y+z\leq r-1$ or $x+y\geq r$, then trivially, respectively, $\mathbf{Pr_0}(A_{I'}|I=I')=0$ or $\mathbf{Pr_1}(A_{I'}|I=I')=1$.
			Analogously, if
			 $x+y+z\geq 2r-1$, then %$x+y\geq r$ or $z\geq r$, and thus, 
			 $\mathbf{Pr_1}(A_{I'}|I=I')=1$, as any $\{a_s,a_t\}$-colouring of $E'$ must
			 repeat at least $r$ times the colour $a_s$ or the colour $a_t$ then.
			 In all these cases, (\ref{P0P1ineq7}) holds. We may therefore assume that 
			 \begin{equation}\label{xyz-ineq1}
			 x+y+z\geq r, \hspace{1cm} x+y\leq r-1, \hspace{1cm} %{\rm and}\hspace{1cm} 
			 x+y+z\leq 2r-2.
			 \end{equation}
Note also that if $y+z\leq r-1$, then $d_{E'_{a_s}}(v)\leq  r-1$ for every admissible $\{a_s,a_t\}$-colouring of $E'$ (from $\mathcal{L}$ or $\mathcal{L}'$).
Thus, $\Omega$ and $\Omega'$ may only include colourings with $d_{E'_{a_t}}(v)\geq  r$, hence by (\ref{xDef}) -- (\ref{zDef'}),  	
$|\Omega'|\geq |\Omega|$,
			since $x+y\geq x$. Similarly, if $x+z\leq r-1$, then we shall always have $d_{E'_{a_t}}(v)\leq  r-1$ while picking colours ($a_s$ or $a_t$) from $\mathcal{L}$. On the other hand, each $\{a_s,a_t\}$-colouring of $E'$ from $\mathcal{L}$ which results in 
$d_{E_{a_s}}(v)\geq  r$ naturally maps bijectively (by mutually interchanging the choices of colours $a_s$ and $a_t$ from the lists corresponding to $z$, and changing from $a_s$ to $a_t$ the colours of all edges corresponding to $y$) to an $\{a_s,a_t\}$-colouring of $E'$ 
from $\mathcal{L}'$ resulting in $d_{E_{a_t}}(v)\geq  r$. 
This again implies that 
$|\Omega'|\geq |\Omega|$. 
			Therefore, we may assume that 	
			 \begin{equation}\label{xyz-ineq1B}
			y+z\geq r \hspace{1cm}{\rm and} \hspace{1cm}x+z\geq r.
			\end{equation}
			
			Note that by the third inequality in~(\ref{xyz-ineq1}), 
			it is not possible that  $d_{E'_{a_s}}(v)\geq  r$ and at the same time $d_{E'_{a_t}}(v)\geq  r$.
 			
			Moreover, by~(\ref{xyz-ineq1B}),  the number of colourings in $\Omega$ with 
			$d_{E'_{a_t}}(v)\geq r$ 
			equals: 
			\begin{equation}\label{Omega1}
			 \binom{z}{r-x}+\binom{z}{r-x+1}+\dotsm+\binom{z}{z}, 
			 \end{equation}
			while in the case of $d_{E'_{a_s}}(v)\geq r$,  
			exactly:
			\begin{equation}\label{Omega2}
			 \binom{z}{r-y}+\binom{z}{r-y+1}+\dotsm+\binom{z}{z}. 
			 \end{equation}		
			
			Analogously, by~(\ref{xyz-ineq1}),  the number of colourings in $\Omega'$ with 
			$d_{E'_{a_t}}(v)\geq r$ 
			equals:
			\begin{equation}\label{Omega'1}
			\binom{z}{r-x-y}+\binom{z}{r-x-y+1}+\dotsm+\binom{z}{z} 
			 \end{equation}
			while in the case of $d_{E'_{a_s}}(v)\geq r$, 
			%exactly (which we interpret as $0$ for $z<r$):
			\begin{equation}\label{Omega'2}
			 \binom{z}{r}+\binom{z}{r+1}+\dotsm+\binom{z}{z} 
			\end{equation}	
				(which we interpret as $0$ if $z\leq r-1$, abusing slightly the use of $\binom{z}{z}$ in~(\ref{Omega'2}) above).
			
			Thus, by (\ref{Omega1}),  (\ref{Omega2}), (\ref{Omega'1}) and (\ref{Omega'2}), in order to show that $|\Omega'|\geq |\Omega|$ it suffices to prove that:
			\begin{equation}	\label{lemma_binom}
				\binom{z}{r-y}+\dotsm+\binom{z}{r-1} \leq \binom{z}{r-x-y}+\dotsm+\binom{z}{r-x-1}
			\end{equation}
			(where some summands on the left-hand side might be equal $0$, if $z\leq r-2$).
			Note that the sums on both sides of the inequality above include $y$ summands each.
			Moreover, by~(\ref{xyz-ineq1}),
			$\frac{(r-y)+(r-1)}{2}-\frac{z}{2}>\frac{x}{2}$, and hence, 
			$|\frac{(r-x-y)+(r-x-1)}{2}-\frac{z}{2}|=|\frac{(r-y)+(r-1)}{2}-\frac{z}{2}-x| < |\frac{(r-y)+(r-1)}{2}-\frac{z}{2}|$.
			  %$\min\{r-y,z-(r-1)\} = z-r+1$.
			%On the other hand, by~(\ref{xyz-ineq1}),  $z-r+1\leq \min\{r-x-y,z-(r-x-1)\}$. 
			Thus, inequality~(\ref{lemma_binom}) follows by Observation~\ref{ObsPascal}.
		\end{proof}

	\end{lemma}

	\begin{theorem}\label{Th:UnifiedVectors}
	For every fixed $\varepsilon\in(0,0.9]$ and an integer $\ell\geq 2$,
	for every $\varepsilon$-excessive tolerance vector
	$\Lambda = (\alpha_1, \alpha_2, \ldots, \alpha_\ell)$ 
	and any graph $G$ with a $\Lambda$-list assignment $L$
	together with the corresponding
	$\Lambda$-tolerance function $\alpha$,
	if $\delta(G)\ge \lceil 109\ell\varepsilon^{-2}\ln(\ell\varepsilon^{-1})\rceil$,
	then there is an $\alpha$-majority $L$-colouring of $G$.
	
	\begin{proof}
	Let us fix $\varepsilon\in(0,0.9]$ and an integer $\ell\geq 2$. Set
	\begin{equation}\label{delta_0_def}
	\delta_1 = \left\lceil 109\ell\varepsilon^{-2}\ln(\ell\varepsilon^{-1})\right\rceil.
	\end{equation}
	Analogously as argued within Observation~\ref{Obs:splitting} on vertex 
	splitting, we may assume that $G=(V,E)$ is a graph with minimum degree 
	$\delta\geq \delta_1$ and maximum degree  $\Delta\leq 2\delta_1-1$.	
	Fix any $\varepsilon$-excessive tolerance vector
	$\Lambda = (\alpha_1, \alpha_2, \ldots, \alpha_\ell)$ and a
	$\Lambda$-list assignment $L:E\rightarrow 2^C$ together with the corresponding
	$\Lambda$-tolerance function $\alpha:C\rightarrow (0,1)$.
	
	For the sake of optimisation and possibly more importantly 
	in order to avoid dependence of $\delta_1$ on $\min(\alpha)$,
	we shall choose parameters within our random approach more carefully than in the proof of Theorem~\ref{Th:MainGeneral}. Let us first define an auxiliary function $h:(0,1)\to (0.5,1)$,
	\begin{equation}\label{h_def}
	h(x)=\frac{\sqrt{1+8x}-1}{4x}.
	\end{equation}
	It is straightforward to verify that $h$ is decreasing in $(0,1)$. Thus, 
	\begin{equation}\label{lambda_def}
	\lambda:=h(0.9)\leq  h(\varepsilon).
	\end{equation}
	Let us set 
	\begin{equation}\label{Mu}
	\mu=\frac{\lambda^2\varepsilon^2}{3\ell}.
	\end{equation}
	In our randomised construction we shall disregard colours with very small tolerance. Without loss of generality we may assume these correspond to entries at the end of the tolerance vector $\Lambda$. More precisely, let us assume that for some $\ell'\leq \ell$:
	\begin{equation}\label{ElPriFeature}
	\alpha_i\geq 6\mu~~~{\rm if~and~only~if}~~~i\in [\ell'].
	\end{equation}
	Note that by~(\ref{Mu}), the sum of small tolerances in $\Lambda$ is limited: 
	\begin{equation}\label{small_tol_sum}
	\sum_{i\in[\ell]\smallsetminus[\ell']}\alpha_i < \ell\cdot 6\mu \leq 2\lambda^2\varepsilon^2.
	\end{equation}
	Hence, by~(\ref{small_tol_sum}), (\ref{lambda_def}), (\ref{h_def})  and~(\ref{AlphaiEps})
	\begin{equation}\label{EpsilonPri}
	\varepsilon':=\varepsilon - 2\lambda^2\varepsilon^2\geq \varepsilon (1-2\lambda^2\cdot0.9)>0,
	\end{equation}
	\begin{equation}\label{AlphaiEpsPrime}
	\sum_{i\in[\ell']}\alpha_i \geq 1+\varepsilon'.
	\end{equation}
	Note also that by~(\ref{EpsilonPri}), (\ref{lambda_def}) and~(\ref{h_def}),
	\begin{eqnarray}
	\varepsilon' &=& \varepsilon - 2\lambda^2\varepsilon^2 
	\geq \varepsilon - 2(h(\varepsilon))^2\varepsilon^2 
	= \varepsilon - \frac{2+8\varepsilon-2\sqrt{1+8\varepsilon}}{8} \nonumber\\
	&=&  \frac{\sqrt{1+8\varepsilon}-1}{4\varepsilon}\cdot\varepsilon
	= h(\varepsilon)\cdot\varepsilon \geq \lambda\varepsilon.
	\label{Lambda-Epsilon}
	\end{eqnarray}
	For every $i\in[\ell']$, let now $\beta_i$ be a real number such that 
	\begin{equation}\label{AlphaiBetai}
	\alpha_i=\beta_i+\sqrt{3\mu}\sqrt{\beta_i},
	\end{equation}
note that 
\begin{equation}\label{Betai_ineq}
0<\beta_i<\alpha_i.
\end{equation}
Moreover, note that by~(\ref{ElPriFeature}), for $\beta_i=3\mu$, we would get that $\beta_i+\sqrt{3\mu}\sqrt{\beta_i}\leq \alpha_i$ for $i\in[\ell']$. Therefore, as $x+\sqrt{3\mu}\sqrt{x}$ is an increasing function, we must have that 
\begin{equation}\label{Betai3Mu}
\beta_i\geq 3\mu
\end{equation}
for  $i\in[\ell']$.
Set 
\begin{equation}\label{BDef}
B=\sum_{i\in [\ell']} \beta_i.
\end{equation} 
By~(\ref{AlphaiEpsPrime}), (\ref{AlphaiBetai}), concavity of $\sqrt{x}$, (\ref{Mu}), 
(\ref{EpsilonPri}) and~(\ref{Lambda-Epsilon}),
\begin{eqnarray}\label{B_Bound}
1+\varepsilon' &\leq& \sum_{i\in[\ell']}\alpha_i 
= \sum_{i\in[\ell']}\beta_i+\sqrt{3\mu}\sum_{i\in[\ell']}\sqrt{\beta_i} \nonumber\\
&=& B + \ell' \sqrt{3\mu}\frac{1}{\ell'}\sum_{i\in[\ell']}\sqrt{\beta_i} 
\leq B + \ell' \sqrt{3\mu}\sqrt{\frac{\sum_{i\in[\ell']}\beta_i}{\ell'}} \nonumber\\
&=& B+\sqrt{\ell'}\frac{\lambda\varepsilon}{\sqrt{\ell}}\sqrt{B}
\leq B + \lambda \varepsilon\sqrt{B}
\leq B + \varepsilon'\sqrt{B}. \nonumber
\end{eqnarray}
Consequently, 
\begin{equation}\label{Batleast1}
B\geq 1.
\end{equation}
 For every $i\in[\ell']$, we set:
\begin{equation}\label{PiDef}
p_i:=\frac{\beta_i}{B}.
\end{equation}
By~(\ref{Betai_ineq}), (\ref{Batleast1}), (\ref{BDef}) and~(\ref{PiDef}),
\begin{eqnarray}
&0<p_i<\alpha_i,&\nonumber\\
&\sum_{i\in [\ell']}p_i=1.&\nonumber 
\end{eqnarray}
Let $X_e$ be a random variable accountable for prescribing a colour to any given edge $e\in E$. Suppose $L(e)=(c_1,c_2,\ldots,c_\ell)$. Then,  for every $c\in C$, the variable $X_e$ takes value $c$ with probability $p_i$ if $c=c_i$ and $i\in[\ell']$ (and $0$ otherwise).
Let $C'\subseteq C$ denote the colours with non-zero probability of being chosen.  For each $e\in E$ and $c\in C'$, we further define the binary random variable $Y_{e,c}$ taking value $1$ if $X_e=c$, and $0$ otherwise. 
	Denote by $L_v$ the set of all colours from $C'$ appearing in lists incident with a vertex $v$, i.e., $L_v=C'\cap(\bigcup_{e\in E(v)}L(e))$. For every $c\in L_v$ we finally again define a random variable:
	$$Z_{v,c}=\sum_{e\in E(v)}Y_{e,c},$$
	representing the number of edges incident with $v$ which are coloured $c$ within our random edge colouring.
	Note that all random variables $Y_{e,c}$ in the sum defining such $Z_{v,c}$ are independent. 
	For $v\in V$ and $c\in L_v$ let us denote events:
	\begin{eqnarray}
	&A_{v,c}:&~~Z_{v,c}>\alpha(c)d(v). \nonumber\\
	&B_v:&~~{\rm therere~exist~} c'\in L_v {\rm ~for~which~} Z_{v,c'}>\alpha(c')d(v). \nonumber
	\end{eqnarray}
	We first want to bound the probability of $B_v$ from the above. 
	%wish to provide an upper bound for $\mathbf{Pr}(B_v)$. 
	By Observation~\ref{The_same_lists}, we may thus assume that all edges $e\in E(v)$ have prescribed exactly the same list, say $L(e)=(c_1,c_2,\ldots,c_\ell)$. Hence, 
	\begin{equation}\label{ListSize}
	|L_v|=\ell'.
	\end{equation} 
	For every $c_i\in L_v$, by~(\ref{PiDef}) and~(\ref{Batleast1}), we have:
	$$\mathbf{E}\left(Z_{v,c_i}\right)=\sum_{e\in E(v)}\mathbf{E}\left(Y_{e,c_i}\right)
	=\sum_{e\in E(v)}p_i = d(v)p_i
	\leq \beta_id(v).$$
	Thus, by~(\ref{AlphaiBetai}) and the Chernoff Bound (applicable due to~(\ref{Betai3Mu})), 
	\begin{eqnarray}
	\mathbf{Pr}\left(A_{v,c_i}\right)
	&=& \mathbf{Pr}\left(Z_{v,c_i}>\alpha_id(v)\right) \nonumber\\
	&=& \mathbf{Pr}\left(Z_{v,c_i}>d(v)\beta_i+d(v)\sqrt{3\mu}\sqrt{\beta_i}\right) 
	\leq e^{-d(v)\mu}. \label{PrAvcBound2} 
	\end{eqnarray}
	Therefore, by~(\ref{PrAvcBound2}) and~(\ref{ListSize}),
	\begin{equation}\label{PrBvBound2}
	\mathbf{Pr}\left(B_v\right)
	\leq  \ell' e^{-d(v)\mu}
	\leq  \ell e^{-d(v)\mu}.
	\end{equation}
	Set 
	 \begin{equation}\label{Gamma1_def}
	 \gamma_1=109\ell\varepsilon^{-2}\ln(\ell\varepsilon^{-1}),
	 \end{equation}
	 hence 
	 \begin{equation}\label{GammaSmall2}
	 \gamma_1\leq \lceil\gamma_1\rceil = \delta_1 \leq \delta \leq d(v)
	 \end{equation}
	 for every vertex $v\in V$.
	 Consider the function 
	 \begin{equation}\label{Gdef}
	 g(x)=2e\ell x e^{-\mu x}
	 \end{equation}
	 and note that
	 $g'(x)=2e\ell e^{-\mu x}(1-\mu x)$, and hence $g$ is decreasing for $x\geq \frac{1}{\mu}$, while by~(\ref{Mu}), (\ref{Gamma1_def}), (\ref{lambda_def}) and~(\ref{h_def}), 
	 $$\gamma_1=\frac{109}{3}\lambda^2\ln(\ell\varepsilon^{-1})\frac{1}{\mu}
	 \geq \frac{109}{3}0.5^2\ln(2\cdot 0.9^{-1})\frac{1}{\mu}>\frac{1}{\mu}.$$
	 Therefore, $g$ is decreasing for $x\geq \gamma_1$. Hence, for every such $x$, by~(\ref{Gdef}), (\ref{Mu}),
	  (\ref{lambda_def}) and~(\ref{h_def}),
\begin{eqnarray}
	 g(x)\leq g(\gamma_1) &=& 
	 2e\ell 109\ell\varepsilon^{-2}\ln(\ell\varepsilon^{-1}) e^{- \frac{109\lambda^2}{3}\ln(\ell\varepsilon^{-1})}\nonumber\\
	 &=& 218e (\ell\varepsilon^{-1})^{2-\frac{109\lambda^2}{3}} \ln(\ell\varepsilon^{-1})  \nonumber\\
	 &\leq& 218e (2\cdot 0.9^{-1})^{2-\frac{109\lambda^2}{3}} \ln(2\cdot 0.9^{-1})    < 1, \label{gxIneq}
\end{eqnarray}
where the last weak inequality above follows by the fact that $g_1(x)=x^{2-\frac{109\lambda^2}{3}}\ln x$ is decreasing for $x$ at least $\exp(-(2-\frac{109\lambda^2}{3})^{-1})<2\cdot 0.9^{-1}$. 
By~(\ref{PrBvBound2}), (\ref{Gdef}), (\ref{gxIneq}) and~(\ref{GammaSmall2}),
\begin{equation}\label{PrBvEstimate2}
\mathbf{Pr}(B_v) \leq \frac{g(d(v))}{2e\cdot d(v)} 
< \frac{1}{2e\cdot d(v)}  
\leq \frac{1}{2e\cdot \delta_1}.  
\end{equation}
Hence, for $p=(2e\cdot\delta_1)^{-1}$, by~(\ref{PrBvEstimate2}), $\mathbf{Pr}(B_v)<p$ for every $v\in V$. Since every event $B_v$ is mutually independent of all but at most $D=d(v)\leq \Delta \leq 2\delta_1-1$ other events $B_{v'}$, by the Lov\'asz Local Lemma, there is a choice of colours from the lists so that none of the events $B_v$, $v\in V$, holds.	

	\end{proof}
	
	\end{theorem}

Note that Theorem~\ref{Th:UnifiedVectors} applied to $1/k$-majority edge colourings from arbitrary lists of length $k+1$ yields a lower bound of order $k^3\log k$ for $\delta$. This matches the bound which could be derived from~\cite{majority23} in such a setting, although the result in Theorem~\ref{Th:UnifiedVectors} concerns a significantly more capacious  concept of majority colourings and a much wider spectrum of list assignments with diversified colour tolerances.

	\section{Concluding remarks\label{Concluding_section}}

There are several directions towards which we 
may further develop %extend 
our research and results.
%{\color{blue} TO JEDYNY FRAGMENT POZOSTALY DO WYMYSLENIA, ZMIANY: In particular, following e.g. an argument from~\cite{list_maj}, based on compactness, one may extend selected theorems from the current paper towards an infinite setting.} 
In particular, using rather standard arguments based on compactness, one may extend selected theorems from the current paper towards an infinite setting, see e.g.~\cite{list_maj} for %examples of 
exemplary instances of such reasonings.
We omit details here, %referring an interested reader to~\cite{list_maj}, 
and only just state one such possible extension. %of our results. % an exemplary extension.
	\begin{theorem}\label{thm_galvin_infinite}
		For every integer $k\geq 2$, each finite and infinite graph $G$ with minimum degree $\delta \geq 2k^2-2k$ has a $1/k$-majority edge colouring from lists of size $k+1$.
	\end{theorem}
		Note that the special case of $k=2$ of this theorem confirms Conjecture~\ref{list_maj_conj_k2} in full length.

	In order to simplify our further discussion though, let us focus our attention back on finite graphs.
	Note that for any $\varepsilon$-excessive tolerance vector $\Lambda = (\alpha_1, \alpha_2, \ldots, \alpha_\ell)$ and a $\Lambda$-list assignment $L$ of a graph $G$  together with the corresponding
	$\Lambda$-tolerance function $\alpha$ such that $\min(\alpha)\geq a$,
	we have $\ell\leq (1+\varepsilon)/a$. 	Consequently, the lower bound for $\delta(G)$ guaranteeing 
	the existence of an  $\alpha$-majority $L$-colouring of $G$ stemming from Theorem~\ref{Th:UnifiedVectors}, and expressed in terms of $\varepsilon$ and $a$ (through substituting $\ell$ with $(1+\varepsilon)/a$) is better with respect to the multiplicative constant than the one in 
 Theorem~\ref{Th:MainGeneral} (for $\varepsilon\in(0,0.9]$). 
 
 In fact the mentioned multiplicative constants in Theorems~\ref{Th:MainGeneral} and~\ref{Th:UnifiedVectors} can be significantly reduced for smaller admissible values of $\varepsilon$, e.g. $\varepsilon\leq 0.1$, which either way seems a more typical setting for our problems.
	Improvements of the bounds in these theorems can also be obtained in the environment of 
	regular graphs. This subclass of graphs admits also a strong improvement of Theorem~\ref{thm_galvin}, which implies confirmation of Conjecture~\ref{EdgeMajconj} in the following special case.	
			\begin{observation}\label{SpecialRegularListK}
		Let $k\geq 2$ be an integer. If $G$ is a $d$-regular graph with $d\geq k^2-k$ such that $\lfloor d/k \rfloor$ is even. Then $G$ has a $1/k$-majority edge colouring from any lists of size $k+1$.
	\end{observation}
	\begin{proof}
	This observation follows directly by the proof of Theorem~\ref{thm_galvin}.
		Since $G$ has only vertices of degree $d$ and  $\lfloor d/k \rfloor$ is even,
		only constraints following from Case~\ref{Case1} 
		of the mentioned proof must be respected, and these boil down to the inequality $t\geq \frac{k-1}{2}$, where $\lfloor d/k\rfloor =2t$,
		which is fulfilled for $d\geq k^2-k$. 
	\end{proof}
This observation supports a direct strengthening of Conjecture~\ref{EdgeMajconj} towards the list setting, which we dare to pose below.
	\begin{conjecture}\label{EdgeMajconjListK}
		For every integer $k\geq 2$, if a graph $G$ has minimum degree $\delta(G)\geq k^2$, then $G$ 
		has a $1/k$-majority edge colouring from any lists of size $k+1$.
	\end{conjecture}
	In a general setting the corresponding lower bounds for $\delta(G)$ in  Theorems~\ref{Th:MainGeneral} and~\ref{tw14} are of the forms $O(a^{-1}\varepsilon^{-2}\ln(a^{-1}\varepsilon^{-1}))$,  $O(\ell^2\varepsilon^{-2})$, while in the case of uniform vector of tolerances, i.e. within Theorem~\ref{Th:UnifiedVectors}, of the form $O(\ell\varepsilon^{-2}\ln(\ell\varepsilon^{-1}))$. We however expect the following should hold true.
	\begin{conjecture}\label{Con:MainGeneralL}
		There is a function $\delta_2:\mathbb{N}\times (0,+\infty)\ni (\ell,\varepsilon)\to\mathbb{R}$ such that $\delta_2=O(\ell\varepsilon^{-1})$ and
	for any fixed $(\ell,\varepsilon)\in \mathbb{N}\times (0,+\infty)$ 
	and every graph $G$ 
	with the minimum degree $\delta(G)\ge \delta_2(\ell,\varepsilon)$, any $\varepsilon$-excessive list assignment $L$ of $G$ associated with arbitrary tolerance function $\alpha$ and containing exclusively list of lengths at most $\ell$, i.e. with $|L(e)|\leq\ell$ for every $e\in E(G)$,
there is an $\alpha$-majority $L$-colouring of $G$.		
\end{conjecture}
		\begin{conjecture}\label{Con:MainGeneralA}
		There is a function $\delta_3:(0,1)\times (0,+\infty)\ni (a,\varepsilon)\to\mathbb{R}$ such that $\delta_3=O(a^{-1}\varepsilon^{-1})$ and
	for any fixed $(a,\varepsilon)\in (0,1)\times (0,+\infty)$ 
	and every graph $G$ 
	with the minimum degree $\delta(G)\ge \delta_3(a,\varepsilon)$, any $\varepsilon$-excessive list assignment $L$ of $G$ 
	and every tolerance function $\alpha$ of $L$ 
	with $\min(\alpha)\geq a$, there is an $\alpha$-majority $L$-colouring of $G$.
	\end{conjecture}
If not in general, it would be interesting to prove these conjectures for uniform vectors of tolerances or for their direct non-list counterparts, which we formalise separately below.

Let $G=(V,E)$ be a graph and $C$ be a set of colours.  
Consider a tolerance function 
$\alpha:C\rightarrow (0,1)$. As before, we denote $\min(\alpha):=\min_{c'\in C}\alpha(c')$.
We call $\alpha$ an  \emph{$\varepsilon$-excessive} tolerance function if 
	$\sum_{c\in C} \alpha(c) \ge 1+\varepsilon$ 
	for some $\varepsilon>0$. 
	A colouring $\omega:E\to C$ is said to be an 
		\emph{$\alpha$-majority $C$-colouring}
		of $G$ if for each vertex $v\in V$ and any colour $c\in C$,
		$d_{E_{c}}(v)\le \alpha(c)d(v)$.
	\begin{conjecture}\label{Con:MainGeneralL2}
		There is a function $\delta_4:\mathbb{N}\times (0,+\infty)\ni (\ell,\varepsilon)\to\mathbb{R}$ such that $\delta_4=O(\ell\varepsilon^{-1})$ and
	for any fixed $(\ell,\varepsilon)\in \mathbb{N}\times (0,+\infty)$ 
	and every graph $G$ 
	with the minimum degree $\delta(G)\ge \delta_4(\ell,\varepsilon)$ and any $\varepsilon$-excessive  tolerance function $\alpha:C\to(0,1)$ with $|C|\leq \ell$ there is an $\alpha$-majority $C$-colouring of $G$.		
\end{conjecture}
		\begin{conjecture}\label{Con:MainGeneralA2}
		There is a function $\delta_5:(0,1)\times (0,+\infty)\ni (a,\varepsilon)\to\mathbb{R}$ such that $\delta_5=O(a^{-1}\varepsilon^{-1})$ and
	for any fixed $(a,\varepsilon)\in (0,1)\times (0,+\infty)$ and every graph $G$ 
	with the minimum degree $\delta(G)\ge \delta_5(a,\varepsilon)$ and any $\varepsilon$-excessive  tolerance function $\alpha:C\to(0,1)$ with $\min(\alpha)\geq a$, there is an $\alpha$-majority $C$-colouring of $G$.
	\end{conjecture}
Still the best known results in such a non-list setting with diversified tolerances for distinct colours follow from
Theorem~\ref{tw14} and  Theorem~\ref{Th:UnifiedVectors}, whose usefulness gets more prominent 
in this setting, where it is applicable and yields a better result than Theorem~\ref{Th:MainGeneral}.
Pushing  down the bounds for $\delta(G)$ which stem from these theorems, and possibly proving Conjectures~\ref{Con:MainGeneralL2} or~\ref{Con:MainGeneralA2} seems to require developing new tools and approaches, though.
A minor support of our suspicion that the non-list setting should be easier to handle is included in Corollary~\ref{FrugalCorollary} below, which almost literally extends Observation~\ref{SpecialRegularListK} to all regular graphs (in a non-list setting), thus confirming Conjecture~\ref{EdgeMajconj} in such a case. It however follows directly by results concerning so-called defective colourings.
			An \textit{edge colouring with defect $d$} or a \textit{$d$-frugal edge colouring} of a graph 
			$G$ is a colouring of the edges of $G$ such that each vertex is incident with at most $d$ edges of the same colour. The minimum $k$ such that $G$ has an edge colouring with defect $d$ using $k$ colours is the \textit{$d$-defective chromatic index} of $G$ and is denoted by $\chi_d^\prime(G)$.
	
	\begin{theorem}[\cite{defective}]\label{simple_defective}
		For every $d\geq 1$ and any graph $G$ with maximum degree $\Delta$, $\chi_d^\prime(G) \in \left\{ \lceil \frac{\Delta}{d} \rceil, \lceil \frac{\Delta+1}{d} \rceil \right\} $.
	\end{theorem}
	
	\begin{corollary}\label{FrugalCorollary}
		If $G$ is an $r$-regular graph with $r \geq k^2$, then $G$ has a $1/k$-majority edge colouring with $k+1$ colours.
	\end{corollary}
	\begin{proof}
		Since $G$ is $r$-regular, every vertex allows $\lfloor r/k \rfloor$ incident edges of the same colour. Hence, $1/k$-majority edge colouring of $G$ is equivalent to its edge colouring with defect $\lfloor r/k \rfloor$.
		Let $r = kt + i$, where $0\leq i<k$ (hence, $t\geq k$). Then $\lfloor r/k \rfloor = t$.
		By Theorem \ref{simple_defective}, such a colouring can thus be constructed by means of at most
		\[ \left\lceil \frac{r+1}{t} \right\rceil 
		= \left\lceil \frac{kt+i+1}{t} \right\rceil
		= k + \left\lceil \frac{i+1}{t} \right\rceil
		\leq k + \left\lceil \frac{k}{t} \right\rceil
		= k+1
		\]
		colours.
	\end{proof}

	Let us finally note that the method used in Theorem \ref{thm_galvin} provides surprisingly strong implications within a non-list setting. Suppose we are interested in a kind of equitable edge colourings of graphs, i.e. such colourings within which each colour is assigned to almost the same number of edges incident with any given vertex. More precisely, consider an edge $k$-colouring of a graph $G=(V,E)$, i.e. an assignment $\omega:E\to[k]$. By a \emph{$k$-discrepancy of $G$}, we understand the following minimum over all $k$-colourings $\omega$ of $G$: 
	$$\mathcal{D}_k(G)=\min_{\omega}\max_{v\in V}\max_{1\leq i<j\leq k}|d_{E_i}(v)-d_{E_j}(v)|.$$
	For $k=2$ this represents a classic discrepancy of a hypergraph with vertex set $E$ and edge set  $\{E(v):v\in V\}$, c.f.~\cite{Lovasz-discrepancy,Spencer-discrepancy} or e.g.~\cite{Balogh-discrepancy} for more recent results in a wider setting. Note that a triangle implies that $D_2(G)$ can be at least $2$. It occurs that this value is always attainable for every graph and each $k$.
	\begin{observation}
For every graph $G$ and each positive integer $k$,  $\mathcal{D}_k(G)\leq 2$.

	\begin{proof}
	This observation follows by yet one more application of the vertex splitting technique, 
	similar as utilized within the proofs of Theorems~\ref{thm_galvin}, \ref{thm_galvin_alpha} and
	Observation~\ref{Obs:splitting}. For every given $G$ and $k$, it suffices to split every vertex $v$ of degree $d(v)$ to $s=\lceil d(v)/(2k)\rceil$ copies so that $s-1$ ones of these copies have degree exactly $2k$ and the remaining one has degree at most $2k$.  The resulting $\overline{G}$ is then, exactly as in the proof of Theorems~\ref{thm_galvin}, made use to define the corresponding bipartite graph $H$ of maximum degree at most $k$. We then colour the edges of $H$ properly with $k$ colours. This colouring is then bijectively reflected through $\overline{G}$ to $G$, as in the mentioned proof. As every vertex of degree $k$ in $H$ was incident with exactly one edge in each colour, then the counterpart of such a  vertex, of degree $2k$ in $\overline{G}$ got incident with exactly two edges in each colour. Consequently, every vertex $v$ in $G$, whose all $s=\lceil d(v)/(2k)\rceil$, except possibly one, copies in $\overline{G}$ were of degree $2k$, must be incident with $2(s-1)$, $2(s-1)+1$ or $2(s-1)+2$ edges in each of the $k$ colours. Thus, the result follows.	
			\end{proof}

\end{observation}

We leave open a problem of determining for which graphs $G$ and which integers $k$, we have  $\mathcal{D}_k(G)\leq 1$.

\end{document}